\documentclass[reqno]{amsart}
\input epsf
\usepackage{amssymb}
\usepackage{epsf}

\newtheorem{theorem}{Theorem}

\theoremstyle{definition}
\newtheorem*{definition}{Definition}

\theoremstyle{remark}
\newtheorem{remark}{Remark}

\begin{document}

\title[Theta Hypergeometric Series]
{Theta Hypergeometric Series}

\author{V.P. Spiridonov}
 \address{Bogoliubov Laboratory of Theoretical Physics, Joint
 Institute for Nuclear Research, Dubna, Moscow Region 141980,
 Russia. }

 \thanks{
Work supported in part by the Russian Foundation for Basic
 Research (RFBR) grant no. 00-01-00299.
\\  \indent Published in:
 Proceedings of the NATO ASI {\em Asymptotic Combinatorics with
 Applications to Mathematical Physics} (St. Petersburg, July 9--22, 2001),
Eds. V. Malyshev and A. Vershik, Kluwer, Dordrecht, 2002, pp. 307--327.
}

\begin{abstract}
We formulate general principles of building hypergeometric type
series from the Jacobi theta functions that generalize the plain
and basic hypergeometric series. Single and
multivariable elliptic hypergeometric series are considered
in detail. A characterization theorem for a single variable
totally elliptic hypergeometric series is proved.
\end{abstract}

\maketitle

\tableofcontents

\section{Introduction}

This note is a mostly conceptual work reflecting partially the content
of a lecture on special functions of hypergeometric type associated with
elliptic beta integrals presented by the author at the NATO Advanced Study
Institute ``Asymptotic Combinatorics with Applications to Mathematical
Physics" (St. Petersburg, July 9-22, 2001). It precedes a
forthcoming extended and technically more elaborate review
\cite{spi:memoir}. Here we describe general principles of building
hypergeometric type series associated with the Jacobi theta functions
\cite{mum:tata}. Some other essential results of \cite{spi:memoir}
were briefly presented in \cite{spi:elliptic,spi:special}.
We discuss only single variable and multivariable series built out of
the Jacobi theta functions though some generalizations based upon the
multidimensional Riemann theta functions are possible.
Moreover, the main attention will be paid to such series obeying certain
ellipticity conditions, i.e. to elliptic hypergeometric series.
We start from a description of the Jacobi theta functions properties
\cite{abr-ste:handbook,whi-wat:course}.

Let us take two complex variables $p$ and $q$ lying inside the
unit disk, i.e. $|p|, |q|<1$. The {\em modular parameters}
$\sigma,\, \text{Im} (\sigma ) > 0,$ and $\tau,\, \text{Im} (\tau)>0,$
are introduced through an exponential representation
\begin{equation}\label{mod-par}
p=e^{2\pi i\tau }, \qquad q=e^{2\pi i\sigma}.
\end{equation}
Define the $p$-shifted factorials \cite{gas-rah:basic}
$$
(a;p)_\infty=\prod_{n=0}^\infty(1-ap^n),\qquad
(a;p)_s=\frac{(a;p)_\infty}{(ap^s;p)_\infty}.
$$
For a positive integer $n$ one has
$$ (a;p)_n=(1-a)(1-ap)\cdots(1-ap^{n-1}) $$
and
$$(a;p)_{-n}=\frac{1}{(ap^{-n};p)_n}. $$
It is convenient to use the following shorthand notations
$$
(a_1, \dots, a_k;p)_\infty \equiv (a_1;p)_\infty\cdots(a_k;p)_\infty.
$$

Let us introduce a Jacobi-type theta function
\begin{equation}\label{theta}
\theta(z;p)=(z,pz^{-1};p)_\infty.
\end{equation}
It obeys the following simple transformation properties
\begin{equation}\label{fun-rel}
\theta(pz;p)=\theta(z^{-1};p)=-z^{-1}\theta(z;p).
\end{equation}
One has also $\theta(p^{-1}z;p)=-p^{-1}z\theta(z;p)$.
Evidently, $\theta(z;p)=0$ for $z=p^{-M},\, M\in\mathbb{Z}$,
and $\theta(z;0)=1-z$.

The standard Jacobi's $\theta_1$-function \cite{whi-wat:course}
is expressed through $\theta(z;p)$ as follows
\begin{eqnarray} \nonumber
\theta_1(u;\sigma,\tau) &=&
-i\sum_{n=-\infty}^\infty (-1)^np^{(2n+1)^2/8}q^{(n+1/2)u}
\\ \nonumber
&=& 2\sum_{n=0}^\infty(-1)^n e^{\pi i\tau(n+1/2)^2}\sin\pi(2n+1)\sigma u
\\ \nonumber
&=& 2p^{1/8}\sin\pi \sigma u\:(p,pe^{2\pi i\sigma u},
pe^{-2\pi i\sigma u};p)_\infty  \\
&=& p^{1/8} iq^{-u/2}\: (p;p)_\infty\: \theta(q^{u};p), \quad
u\in\mathbb{C}.
\label{theta1}\end{eqnarray}
We have introduced artificially the second modular parameter $\sigma$
into the definition of $\theta_1$-function---the variable
$u$ will often take integer values and it is convenient to
make an appropriate rescaling from the very beginning.
Note that other Jacobi theta functions $\theta_{2,3,4}(u)$
can be obtained from $\theta_1(u)$ by a simple shift of the
variable $u$ \cite{whi-wat:course}, i.e. these functions
structurally do not differ much from $\theta_1(u)$.

In the following considerations we shall be employing
convenient notations by replacing $\theta_1$-symbol in favor
of the elliptic numbers $[u]$ used in \cite{djkmo:exactly}:
\begin{eqnarray*}
&& [u] \equiv \theta_1(u) \qquad \text{or} \qquad
[u;\sigma,\tau]\equiv \theta_1(u;\sigma,\tau), \\
&& \makebox[4em]{}
[u_0,\ldots,u_k]\equiv \prod_{m=0}^k[u_m].
\end{eqnarray*}
Dependence on $\sigma$ and $\tau$ will be indicated explicitly only if
it is necessary. The function $[u]$ is entire,
odd $[-u]=-[u]$, and doubly quasiperiodic
\begin{eqnarray}\nonumber
&& [u+\sigma^{-1}]     = -[u],   \\
&& [u+\tau\sigma^{-1}] = -e^{-\pi i\tau-2\pi i\sigma u}[u].
\label{quasi}\end{eqnarray}
It is well-known that the theta function $[u]$ can be derived
uniquely (up to a constant factor) from the transformation properties
(\ref{quasi}) and the demand of entireness.

Modular transformations are described by the following $SL(2,\mathbb{Z})$
group action upon the modular parameters $\sigma$ and $\tau$
\begin{equation}\label{sl2}
\tau\to \frac{a\tau+b}{c\tau+d},\qquad \sigma\to \frac{\sigma}{c\tau+d},
\end{equation}
where $a,b,c,d \in\mathbb{Z}$ and $ad-bc=1$. This group is
generated by two simple transformations: $\tau\to\tau+1$,
$\sigma\to\sigma,$ and $\tau\to -\tau^{-1},$ $\sigma\to\sigma\tau^{-1}$.
In these two cases one has
\begin{eqnarray}\label{modular-tr1}
&& [u;\sigma,\tau+1] = e^{\pi i/4} [u;\sigma,\tau],\qquad
\\  \label{modular-tr2}
&& [u;\sigma/{\tau},-1/\tau]
= i(-i\tau)^{1/2}e^{\pi i\sigma^2u^2/\tau} [u;\sigma,\tau],
\end{eqnarray}
where the square root sign of $(-i\tau)^{1/2}$ is fixed from the
condition that the real part of this expression is positive.

\section{Theta hypergeometric series $_rE_s$ and $_rG_s$}

Now we are going to introduce formal power series $_rE_s$ and $_rG_s$
built out from the Jacobi theta functions (we will not consider their
convergence properties here). They generalize the
plain hypergeometric series $_rF_s$ and basic hypergeometric series
$_r\Phi_s$ together with their bilateral partners $_rH_s$ and $_r\Psi_s$.
The definitions given below follow the general
spirit of qualitative (but constructive) definitions of the
plain and basic hypergeometric series going back to Pochhammer and
Horn \cite{aar:special,gas-rah:basic,ggr:general}. In order to
define theta hypergeometric series we use as a key property
the quasiperiodicity of Jacobi theta functions (\ref{quasi}).

\begin{definition}
The series $\sum_{n\in\mathbb{N}} c_n$ and $\sum_{n\in\mathbb{Z}} c_n$
are called {\em theta hypergeometric series of elliptic type}
if the function $h(n)=c_{n+1}/{c_n}$ is a meromorphic doubly quasiperiodic
function of $n$ considered as a complex variable. More precisely, for
$x\in \mathbb{C}$ the function $h(x)$ should obey the following properties:
\begin{equation}\label{h(n)-period}
h(x+\sigma^{-1})=ah(x),\qquad
h(x+\tau\sigma^{-1})=be^{2\pi i\sigma \gamma x}h(x),
\end{equation}
where $\sigma^{-1}, \tau\sigma^{-1}$ are quasiperiods of the theta function
$[u]$ (\ref{quasi}) and $a, b, \gamma$ are some complex numbers.
\end{definition}

\begin{theorem}
Let a meromorphic function $h(n)$ satisfies the properties (\ref{h(n)-period}).
Then it has the following general form in terms of the Jacobi theta functions
\begin{equation}
h(n)=\frac{[n+u_1,\ldots, n+u_r]}{[n+v_1,\ldots, n+v_s]}\, q^{\beta n} y,
\label{h(n)-quasi}\end{equation}
where $r, s$ are arbitrary non-negative integers and
$u_1,\ldots, u_r, v_1,\ldots,v_s, \beta, y$ are arbitrary complex
parameters restricted by the condition of non-singularity of $h(n)$
and related to the quasiperiodicity multipliers $a, b, \gamma$ as follows:
\begin{eqnarray}\nonumber
&& a = (-1)^{r-s} e^{2\pi i\beta}, \qquad \gamma= s-r, \\
&& b=(-1)^{r-s}e^{\pi i\tau (s-r+2\beta)}
e^{2\pi i\sigma(\sum_{m=1}^sv_m-\sum_{m=1}^ru_m)}.
\label{multipliers}\end{eqnarray}
\end{theorem}
\begin{proof}
Let us tile the complex plane of $n$ by parallelograms whose edges
are formed by the theta function quasiperiods $\sigma^{-1}$ and $\tau\sigma^{-1}$.
Because the multipliers in the quasiperiodicity conditions (\ref{h(n)-period})
are entire functions of $n$ without zeros, the meromorphic function
$h(n)$ has the same finite number of zeros and poles in each
parallelogram on the plane.
Let us denote as $-u_1,\ldots, -u_r$ the zeros of $h(n)$ in one
of such parallelograms and as $-v_1,\ldots, -v_s$ its poles.
For simplicity we assume that these zeros and poles are simple---
for creating multiple zeros or poles it is sufficient to set some
of the numbers $u_m$ or $v_m$ to be equal to each other.

Let us represent the ratio of theta hypergeometric series
coefficients as follows: $c_{n+1}/c_n=h(n)g(n)$, where $h(n)$
has the form (\ref{h(n)-quasi}) with some unfixed parameter $\beta$.
Since all the zeros and poles of $c_{n+1}/c_n$ are sitting in
$h(n)$, the function $g(n)$ must be an entire function without
zeros satisfying the constraints $g(n+\sigma^{-1})=a'g(n),$
$g(n+\tau\sigma^{-1})=b'e^{2\pi i\sigma \gamma' n}g(n)$ for some
complex numbers $a', b', \gamma'$. However, the only function satisfying
such demands is the exponential $q^{\beta' n}$, where $\beta'$
is a free parameter. Since a factor of such type is already present in
(\ref{h(n)-quasi}) we may set $g(n)=1$ and this proves that the
most general $c_{n+1}/c_n$ has the form (\ref{h(n)-quasi}) (note
that $y$ is just an arbitrary proportionality constant). Direct
application of the properties (\ref{quasi}) yields the connection
between multipliers $a, b, \gamma$ and the parameters $u_m,\,
m=1,\ldots,r,$ $v_k,\, k=1,\ldots,s,$ $\beta$ as stated in
(\ref{multipliers}).
\end{proof}

Resolving the first order recurrence relation for the series coefficients
$c_n$ and normalizing $c_0=1$ we get the following explicit ``additive"
expression for the general theta hypergeometric series of elliptic type
\begin{equation}
\sum_{n\in\mathbb{N}\, \text{or}\, \mathbb{Z}}
\frac{[u_1,\ldots,u_{r}]_n}
{[v_1,\ldots,v_s]_n}\,q^{\beta n(n-1)/2} y^n,
\label{_rE_s-1}\end{equation}
where we used the elliptic shifted factorials defined
for $n\in\mathbb{N}$ as follows:
\begin{eqnarray}
&& [u_1,\ldots,u_k]_{\pm n}=\prod_{m=1}^k[u_m]_{\pm n},
\nonumber \\
&& [u]_n=[u][u+1]\cdots[u+n-1],\qquad [u]_{-n}=\frac{1}{[u-n]_n}.
\label{ell-shift}\end{eqnarray}

In order to simplify the trigonometric degeneration limit
$\text{Im}(\tau)\to+\infty$ (or $p\to 0$) in the series (\ref{_rE_s-1})
we renormalize $y$ and introduce as a main series argument another variable $z$:
$$
y \equiv (ip^{1/8})^{s-r}q^{(u_1+\ldots+u_r-v_1-\ldots-v_s)/2}z.
$$
Let us replace the parameter $\beta$ by another parameter $\alpha$ as well
through the relation $\beta\equiv \alpha + (r-s)/2$.
Then, we can rewrite the function (\ref{h(n)-quasi}) in the following
``multiplicative" form using the functions $\theta(tq^n;p)$:
\begin{equation}
h(n)=\frac{\theta(t_1q^n,\ldots,t_{r}q^n;p)}
{\theta(w_1q^n,\ldots, w_{s}q^n;p)}\, q^{\alpha n}z,
\label{h(n)-2}\end{equation}
where $t_m=q^{u_m},\, m=1,\ldots, r,$  $w_k=q^{v_k},\, k=1,\ldots, s,$ and
the following shorthand notations are employed:
$$
\theta(t_1,\ldots,t_k;p)=\prod_{m=1}^k\theta(t_m;p).
$$

Now we are in a position to introduce the unilateral theta hypergeometric
series $_rE_s$. In its definition we follow the standard plain and basic
hypergeometric series conventions. Namely, in the expression (\ref{h(n)-2})
we replace $s$ by $s+1$ and set $u_r\equiv u_0$ and $v_{s+1}=1$. This does
not restrict generality of consideration since one can remove such a
constraint by fixing one of the numerator parameters $u_m$ to be equal to 1.
Then we fix
\begin{eqnarray}\nonumber
\lefteqn{ _rE_s\left({t_0,\ldots, t_{r-1}\atop w_1,\ldots,w_s};q,p;\alpha,
z\right) } && \\ && \makebox[4em]{}
= \sum_{n=0}^\infty \frac{\theta(t_0,t_1,\ldots,t_{r-1};p;q)_n}
{\theta(q,w_1,\ldots,w_s;p;q)_n}\, q^{\alpha n(n-1)/2} z^n,
\label{_rE_s-2}\end{eqnarray}
where we have introduced new notations for the elliptic shifted factorials
$$
\theta(t;p;q)_n=\prod_{m=0}^{n-1}\theta(tq^m;p)
$$
and
$$
\theta(t_0,\ldots,t_k;p;q)_n=\prod_{m=0}^k\theta(t_m;p;q)_n.
$$
We draw attention to the particular ordering of $q,p$
used in the notations for theta hypergeometric series (in the
previous papers we were ordering $q$ after $p$ which does not match
with the ordering in terms of modular parameters $\sigma, \tau$
in $[u;\sigma,\tau]$).

An important fact is that theta hypergeometric series do not admit confluence
limits. Indeed, because of
the quasiperiodicity of theta functions the limits of parameters
$t_m, w_m\to 0$ or $t_m, w_m \to\infty$ are not well defined and it is not possible
to pass in this way from $_rE_s$-series to similar series with smaller
values of indices $r$ and $s$.

For the bilateral theta hypergeometric series we introduce
different notations:
\begin{eqnarray}\nonumber
\lefteqn{ _rG_s\left({t_1,\ldots, t_{r}\atop w_1,\ldots,w_{s}};q,p;\alpha, z
\right) } && \\ && \makebox[4em]{}
=\sum_{n=-\infty}^\infty\frac{\theta(t_1,\ldots,t_{r};p;q)_n}
{\theta(w_1,\ldots,w_{s};p;q)_n}\, q^{\alpha n(n-1)/2} z^n.
\label{_rG_s-2}\end{eqnarray}
This expression was derived with the help of (\ref{h(n)-2})
without any changes. The elliptic shifted factorials for negative
indices are defined in the following way:
$$
\theta(t;p;q)_{-n}=\frac{1}{\theta(tq^{-n};p;q)_n},\quad n\in\mathbb{N}.
$$
Due to the property $\theta(q;p;q)_{-n}=0$ (or $[1]_{-n}=0$) for $n>0$,
the choice $t_{s+1}=q$ (or $v_{s+1}=1$) in the $_rG_{s+1}$ series leads
to its termination from one side. After denoting  $t_r\equiv t_0$ (or
$u_{r}\equiv u_0$) one gets in this way the general $_rE_s$-series.
Since the bilateral
series are more general than the unilateral ones, it is sufficient to
prove key properties of theta hypergeometric series in the
bilateral case without further specification to the unilateral one.

Consider the limit $\text{Im}(\tau)\to+\infty$ or $p\to 0$.
In a straightforward manner one gets
\begin{eqnarray}\nonumber
\lefteqn{ \lim_{p\to 0} {_rE_s} =
{_r\Phi_s}\left({t_0,t_1,\ldots, t_{r-1}
\atop w_1,\ldots,w_s }; q;\alpha, z\right) } && \\ && \makebox[3em]{}
=\sum_{n=0}^\infty \frac{(t_0,t_1,\ldots,t_{r-1};q)_n}
{(q,w_1,\ldots,w_s;q)_n}\, q^{\alpha n(n-1)/2} z^n.
\label{Phi-2}\end{eqnarray}
This basic hypergeometric series is different from the standard one by
the presence of an additional parameter $\alpha$. The definition of
$_r\Phi_s$ series suggested in \cite{sla:generalized} uses $\alpha=0$.
The definition given in \cite{gas-rah:basic} looks as follows
\begin{equation}
{_r\Phi_s}=\sum_{n=0}^\infty \frac{(t_0,t_1,\ldots,
t_{r-1};q)_n}{(q,w_1,\ldots,w_s;q)_n}
\left((-1)^nq^{n(n-1)/2}\right)^{s+1-r}\, z^n,
\label{Phi-3}\end{equation}
which matches with (\ref{Phi-2}) for $\alpha=s+1-r$ after the
replacement of $z$ by $(-1)^{s+1-r}z$. Actually, the $\alpha=0$
and $\alpha=s+1-r$ choices are related to each other through the
inversion transformation $q\to q^{-1}$ with the subsequent
redefinition of parameters $t_m, w_m, z$.
One of the characterizations of the basic hypergeometric series $\sum_nc_n$
consists in the demand for $c_{n+1}/c_n$ to be a general rational function
of $q^n$ which is satisfied by (\ref{Phi-2}) only for integer $\alpha$.
In order to get in the limit $p\to 0$ the standard $q$-hypergeometric
series we fix $\alpha=0$. It is not clear at the moment whether this choice
is the most ``natural" one or it does not play a fundamental role---this
question can be answered only after the discovery of good applications for
the series $_rE_s$ in pure mathematical or mathematical physics problems.
From the point of view of elliptic beta integrals \cite{spi:elliptic,spi:memoir}
this is the most natural choice indeed.

In the bilateral case we fix $\alpha=0$ as well, so that in the $p\to 0$ limit
the $_rG_s$ series are reduced to the general $_r\Psi_s$-series:
\begin{eqnarray}\nonumber
\lefteqn{{_r\Psi_s}\left({t_1,\ldots, t_{r}
\atop w_1,\ldots,w_s }; q;z\right) } && \\ && \makebox[3em]{}
=\sum_{n=-\infty}^\infty \frac{(t_1,\ldots,t_{r};q)_n}
{(w_1,\ldots,w_s;q)_n}\, z^n.
\label{Psi}\end{eqnarray}

\begin{definition}
The series $_{r+1}E_r$ and $_{r}G_{r}$ are called {\em balanced} if
their parameters satisfy the constraints, in the additive form,
\begin{equation}
u_0+\ldots+u_{r}=1+v_1+\ldots+v_r
\label{balance1}\end{equation}
and
\begin{equation}
u_1+\ldots+u_{r}=v_1+\ldots+v_r
\label{balance-g}\end{equation}
respectively.
In the multiplicative form these restrictions look as follows:
$\prod_{m=0}^r t_m$  $=q\prod_{k=1}^rw_k$ and $\prod_{m=1}^r t_m=$
$\prod_{k=1}^rw_k$ respectively.
\end{definition}

\begin{remark}
In the limit $p\to 0$ the series $_{r+1}E_r$ goes to $_{r+1}\Phi_r$
provided the parameters $u_m$ (or $t_m$), $m=0,\ldots,r,$ and
$v_k$ (or $w_k$), $k=1,\ldots,r,$  {\em remain fixed}.
Then our condition of balancing does not
coincide with the one given in \cite{gas-rah:basic}, where $_{r+1}\Phi_r$
is called balanced provided $q\prod_{m=0}^r t_m=\prod_{k=1}^rw_k$
(simultaneously one usually assumes also that $z=q$, but we drop this
requirement). A discrepancy in these definitions will be resolved
after imposing some additional constraints upon the series parameters
(see the very-well-poisedness condition below).
\end{remark}

\section{Elliptic hypergeometric series}

From the author's point of view the following definition plays a
fundamental role for the whole theory of hypergeometric type series
since it explains origins of some known peculiarities of the plain
and basic hypergeometric series.

\begin{definition}
The series $\sum_{n\in\mathbb{N}} c_n$ and $\sum_{n\in\mathbb{Z}} c_n$
are called {\em elliptic hypergeometric series}
if $h(n)=c_{n+1}/{c_n}$ is an elliptic function of the argument $n$
which is considered as a complex variable, i.e. $h(x)$ is a meromorphic
double periodic function of $x\in\mathbb{C}$.
\end{definition}

\begin{theorem}
Let $\sigma^{-1}$ and $\tau\sigma^{-1}$ be two periods of the
elliptic function $h(x)$, i.e. $h(x+\sigma^{-1})=h(x)$
and $h(x+\tau\sigma^{-1})=h(x)$. Let $r+1$ be the order of the elliptic
function $h(x)$, i.e. the number of its poles (or zeros) in the
parallelogram of periods. Then the unilateral (or bilateral) elliptic
hypergeometric series coincides with the balanced theta hypergeometric
series $_{r+1}E_r$ (or $_{r+1}G_{r+1}$).
\end{theorem}
\begin{proof}
It is well known that any elliptic function $h(x),\, x\in\mathbb{C}$,
of the order $r+1$ with the periods $\sigma^{-1}$ and $\tau\sigma^{-1}$
can be written as a ratio of $\theta_1$-functions as follows
\cite{whi-wat:course}:
\begin{equation}
h(x)=z\prod_{m=0}^r\frac{[x+\alpha_m;\sigma,\tau]}{[x+\beta_m;\sigma,\tau]},
\label{factorization}\end{equation}
where the zeros $\alpha_0,\ldots,\alpha_r$ and the poles
$\beta_0,\ldots,\beta_r$ satisfy the following constraint:
\begin{equation}\label{balance2}
\sum_{m=0}^r\alpha_m=\sum_{m=0}^r\beta_m.
\end{equation}
Now the identification of the unilateral elliptic hypergeometric
series with the balanced $_{r+1}E_r$-series is evident. One just has
to shift $x\to x-\beta_0+1$, set $x\in \mathbb{N}$, denote
$u_m=\alpha_m-\beta_0+1,\, v_m=\beta_m-\beta_0+1$, and resolve
the recurrence relation $c_{n+1}=h(n)c_n$. After this, the condition
(\ref{balance2}) becomes the balancing condition for the $_{r+1}E_r$
series. A similar situation takes place, evidently, in the bilateral
series case $_{r+1}G_{r+1}$, when $\alpha_m$ and $\beta_m$ just coincide
with $u_m$ and $v_m$ respectively.

Note that because of the balancing condition (\ref{balance2})
the function (\ref{factorization}) can be rewritten as a simple
ratio of $\theta(t;p)$-functions:
$$
h(x)=z\prod_{m=0}^r\frac{\theta(t_mq^x;p)}{\theta(w_mq^x;p)},
$$
where $t_m=q^{\alpha_m}$ and $w_m=q^{\beta_m}$.
\end{proof}

\begin{definition}
Theta hypergeometric series of elliptic type are called modular
hypergeometric series if they are invariant with respect to the
$SL(2,\mathbb{Z})$ group action (\ref{sl2}).
\end{definition}

Consider what kind of constraints upon the parameters of $_rE_s$
and $_rG_s$ series one has to impose in order to get the modular
hypergeometric series. Evidently, it is sufficient to establish modularity
of the function $h(n)=c_{n+1}/c_n$. From its explicit form (\ref{h(n)-quasi})
and the transformation laws (\ref{modular-tr1}) and (\ref{modular-tr2})
it is easy to see that in the unilateral case one must have
$$
\sum_{m=0}^{r-1}(x+u_m)^2=(x+1)^2+\sum_{m=0}^s(x+v_m)^2,
$$
which is possible only if a) $s=r-1$, b) the parameters satisfy the
balancing condition (\ref{balance1}), and c) the following constraint
is valid:
\begin{equation}
u_0^2+\ldots+u_{r-1}^2=1+v_1^2+\ldots+v_{r-1}^2.
\label{modular-E}\end{equation}
Under these conditions the $_{r}E_{r-1}$-series becomes modular invariant.
Let us note that modularity of theta hypergeometric series assumes
their ellipticity. The opposite is not correct, but a more strong
demand of ellipticity, to be formulated below, automatically leads
to modular invariance. Modular hypergeometric series represent particular
examples of Jacobi modular functions in the sense of Eichler and Zagier
\cite{eic-zag:theory}.

In the bilateral case one must have $r=s$, the balancing condition
(\ref{balance-g}), and the constraint
\begin{equation}
u_1^2+\ldots+u_{r}^2=v_1^2+\ldots+v_r^2
\label{modular-G}\end{equation}
for the $_{r}G_{r}$-series to be modular invariant.

\begin{definition}
The theta hypergeometric series $_{r+1}E_r$ is called {\em well-poised}
if its parameters satisfy the following constraints
\begin{equation}\label{well-poised-1}
u_0+1=u_1+v_1=\ldots=u_{r}+v_r
\end{equation}
in the additive form or
\begin{equation}\label{well-poised-2}
qt_0=t_1w_1=\ldots=t_{r}w_r
\end{equation}
in the multiplicative form.
Similarly, the series $_{r}G_{r}$ is called well-poised if
$u_1+v_1=\ldots=u_r+v_r$ or $t_1w_1=\ldots=t_rw_r$.
\end{definition}

This definition of well-poised series matches with the
one used in the theory of plain and basic hypergeometric series
\cite{gas-rah:basic}. Note that it does not imply the balancing
condition.

\begin{definition}
The series $_{r+1}E_r$ is called {\em very-well-poised}
if, in addition to the constraints (\ref{well-poised-1})
or (\ref{well-poised-2}), one imposes the restrictions
\begin{eqnarray}\nonumber
&& u_{r-3}=\frac{1}{2}u_0+1,\quad u_{r-2}=\frac{1}{2}u_0+1-\frac{1}{2\sigma}, \\
&& u_{r-1}=\frac{1}{2}u_0+1-\frac{\tau}{2\sigma},
\quad u_{r} = \frac{1}{2}u_0+1+\frac{1+\tau}{2\sigma},
\label{very-well-poised-1}\end{eqnarray}
or, in the multiplicative form,
\begin{eqnarray}\nonumber
&& t_{r-3}=t_0^{1/2}q,\quad t_{r-2}=-t_0^{1/2}q, \\
&& t_{r-1}=t_0^{1/2}qp^{-1/2}, \quad t_{r} =- t_0^{1/2}qp^{1/2}.
\label{very-well-poised-2}\end{eqnarray}
\end{definition}

Let us derive a simplified form of the very-well-poised series.
First, we notice that
$$
\theta(zp^{-1/2};p)=-zp^{-1/2}\theta(zp^{1/2};p)
$$
and
$$
\theta(z,-z,zp^{1/2},-zp^{1/2};p)=\theta(z^2;p).
$$
After application of these relations, one can find that
$$
\frac{\theta(t_{r-3},\ldots,t_{r};p;q)_n}
{\theta(qt_0/t_{r-3},\ldots,qt_0/t_{r};p;q)_n}=
\frac{\theta(t_0q^{2n};p)}{\theta(t_0;p)}\, (-q)^n.
$$
As a result, one gets
\begin{eqnarray}\nonumber
\lefteqn{
_{r+1}E_r\left({t_0,t_1,\ldots, t_{r-4},qt_0^{1/2},-qt_0^{1/2},
qp^{-1/2}t_0^{1/2},-qp^{1/2}t_0^{1/2} \atop
qt_0/t_1,\ldots,qt_0/t_{r-4},t_0^{1/2},-t_0^{1/2},
p^{1/2}t_0^{1/2},-p^{-1/2}t_0^{1/2} };q,p;z\right) } &&
\\ && \makebox[5em]{}
= \sum_{n=0}^\infty \frac{\theta(t_0q^{2n};p)}{\theta(t_0;p)}
\prod_{m=0}^{r-4}\frac{\theta(t_m;p;q)_n}{\theta(qt_0/t_m;p;q)_n}\, (-qz)^n.
\label{vwp-1}\end{eqnarray}

For convenience we introduce separate notations for the very-well-poised
series, since they contain an essentially smaller number of parameters
than the general theta hypergeometric series $_{r+1}E_r$.
For this we replace $z$ by $-z$ and all the parameters $t_m,\, m=0,
\ldots,r-4,$ by $t_0t_m$ (in particular, this replaces $t_0$ by $t_0^2$).
Then we write
\begin{eqnarray} \nonumber
\lefteqn{
_{r+1}E_r(t_0;t_1,\ldots,t_{r-4};q,p;z) }&&
\\ && \makebox[3em]{}
\equiv \sum_{n=0}^\infty \frac{\theta(t_0^2q^{2n};p)}{\theta(t_0^2;p)}
\prod_{m=0}^{r-4}\frac{\theta(t_0t_m;p;q)_n}{\theta(qt_0t_m^{-1};p;q)_n}\,
(qz)^n.  \label{vwp-2}\end{eqnarray}
In terms of the elliptic numbers this series takes the form:
\begin{eqnarray} \nonumber
\lefteqn{
_{r+1}E_r(u_0;u_1,\ldots,u_{r-4};\sigma,\tau;z) }&&
\\ && \makebox[3em]{}
\equiv \sum_{n=0}^\infty \frac{[2u_0+2n]}{[2u_0]}\prod_{m=0}^{r-4}
\frac{[u_0+u_m]_n}{[u_0+1-u_m]_n}\, z^nq^{n(\sum_{m=0}^{r-4}u_m-(r-7)/2)}.
\label{vwp-3}\end{eqnarray}
We use the same symbol $_{r+1}E_r$ in (\ref{vwp-2}) and (\ref{vwp-3})
since these series can be easily distinguished from
the general $_{r+1}E_r$ series by the number of free parameters.

For $p=0$ theta hypergeometric series (\ref{vwp-2}) are reduced to the
very-well-poised basic hypergeometric series:
\begin{eqnarray} \nonumber
\lefteqn{
_{r-1}\Phi_{r-2}(t_0;t_1,\ldots,t_{r-4};q;qz) }&&
\\ \nonumber && \makebox[3em]{}
= \sum_{n=0}^\infty \frac{1-t_0^2q^{2n}}{1-t_0^2}
\prod_{m=0}^{r-4}\frac{(t_0t_m;q)_n}{(qt_0t_m^{-1};q)_n}\, (qz)^n,
\label{vwp-phi}\end{eqnarray}
which are different from the corresponding partners in \cite{gas-rah:basic}
by the replacement of $z$ by $qz$ and $t_m$ by $t_0t_m$
(a standard notation for these series would be $_{r-1}W_{r-2}$
but we are not using it here).

Remind that the balancing condition is not involved into the definition
of the very-well-poised theta hypergeometric series. Imposing
the corresponding constraint
$$
\sum_{m=0}^r(u_0+u_m)=1+\sum_{m=1}^r(u_0+1-u_m)
$$
upon (\ref{vwp-3}) we get
$$
\sum_{m=0}^{r-4} u_m= \frac{r-7}{2}.
$$
In the multiplicative form this condition takes the form
$$
\prod_{m=0}^{r-4}t_m= q^{(r-7)/2}.
$$
But this is precisely
the balancing condition for the very-well-poised series
appearing in the theory of basic hypergeometric series
\cite{gas-rah:basic}. Thus for the very-well-poised series
there is no discrepancy in the definitions of balancing condition
given in \cite{gas-rah:basic} and in this paper. This happens
because the constraints (\ref{very-well-poised-1}) taken separately
are not well defined in the limit $\text{Im} (\tau)\to+\infty$.
Note that for the balanced series an extra factor standing in
(\ref{vwp-3}) to the right of $z^n$ disappears.
Summarizing this consideration we conclude that a very natural
condition of ellipticity of the function $h(n)=c_{n+1}/c_n$
in the theta hypergeometric series provides a substantial
meaning to the (innatural) balancing condition for the standard
basic hypergeometric series.

If we impose balancing condition in the multiplicative form
then there appears an ambiguity.
Indeed, substituting into the condition $\prod_{m=0}^r(t_0t_m)=
q\prod_{m=1}^r(qt_0/t_m)$ the constraints
$t_{r-3}=q,\, t_{r-2}=-q,\, t_{r-1}=qp^{-1/2},\, t_{r} = -qp^{1/2}$
(these are the restrictions (\ref{very-well-poised-2}) after the shift
$t_m\to t_0t_m$) we get $\prod_{m=0}^{r-4}t_m^2=q^{r-7}$ which
yields $\prod_{m=0}^{r-4}t_m=\pm q^{(r-7)/2}$ and it is known that
only the plus sign corresponds to the correct balancing condition
for odd $r$ (the even $r$ cases remain ambiguous, but even $r$
do not appear in known examples of summation formulae of basic
hypergeometric series).

In the same way, the bilateral theta hypergeometric series
$_{r}G_{r}$ are called very-well-poised if the constraints
(\ref{very-well-poised-1}) or (\ref{very-well-poised-2})
are satisfied, where $u_0$ or $t_0=q^{u_0}$ is a free parameter.
Following the unilateral series case, we
replace $z$ by $-z$, shift the parameters $t_m\to t_0t_m,\,
m=0,\ldots, r-4,$ and introduce the following shorthand notations
for the simplified form of these series:
\begin{eqnarray} \nonumber
\lefteqn{
_{r}G_{r}(t_0;t_1,\ldots,t_{r-4};q,p;z) }&&
\\  && \makebox[3em]{}
= \sum_{n=-\infty}^\infty \frac{\theta(t_0^2q^{2n};p)}{\theta(t_0^2;p)}
\prod_{m=1}^{r-4}\frac{\theta(t_0t_m;p;q)_n}{\theta(qt_0t_m^{-1};p;q)_n}\,
(qz)^n  \label{vwp-g2}\end{eqnarray}
or in terms of the elliptic numbers
\begin{eqnarray} \nonumber
\lefteqn{
_{r}G_{r}(u_0;u_1,\ldots,u_{r-4};\sigma,\tau ;z) }&&
\\ && \makebox[3em]{}
= \sum_{n=-\infty}^\infty \frac{[2u_0+2n]}{[2u_0]}\prod_{m=1}^{r-4}
\frac{[u_0+u_m]_n}{[u_0+1-u_m]_n}\, z^nq^{n(\sum_{m=1}^{r-4}u_m-(r-8)/2)}.
\label{vwp-g3}\end{eqnarray}

Repeating considerations for the bilateral series we find
the following compact form for the balancing condition for the
$_rG_r$-series (\ref{vwp-g2}):
$$
\sum_{m=1}^{r-4}u_m=\frac{r-8}{2}, \quad
\text{or} \quad \prod_{m=1}^{r-4}t_m=q^{(r-8)/2}.
$$

Under the constraint $u_{r-3}=u_0$ (or $t_{r-3}=t_0$) the
$_{r+1}G_{r+1}$-series is converted into the $_{r+1}E_r$ series.
A general connection between the very-well-poised series of $E$
and $G$ types looks as follows:
\begin{eqnarray}\nonumber
\lefteqn{ _{r}G_{r}(t_0;t_1,\ldots,t_{r-4};q,p;z)=
{_{r+2}E_{r+1}}(t_0;t_1,\ldots,t_{r-4},qt_0^{-1};q,p;z)} && \\ \label{G/E}
&& \makebox[4em]{}
+\frac{q^{r-7}}{z\prod_{m=1}^{r-4}t_m^2}
\frac{\theta(t_0^{-2}q^2;p)}{\theta(t_0^{-2};p)}
\prod_{m=1}^{r-4}\frac{\theta(t_mt_0^{-1};p)}{\theta(qt_0^{-1}t_m^{-1};p)}
\\  \nonumber && \makebox[4em]{}
\times {_{r+2}E_{r+1}}\left(\frac{q}{t_0};t_1,\ldots,
t_{r-4},t_0;q,p; \frac{q^{r-8}}{z\prod_{m=1}^{r-4}t_m^2}\right).
\end{eqnarray}

\begin{remark}
Within our classifications, an elliptic generalization of basic
hypergeometric series $_{r+1}\Phi_r$ introduced by Frenkel and Turaev
in \cite{fre-tur:elliptic} coincides with the very-well-poised
balanced theta hypergeometric series $_{r+1}E_r$
of the unit argument $z=1$. Such series have their origins in
elliptic solutions of the Yang-Baxter equation \cite{bax:exactly,abf:eight,
djkmo:exactly,fre-tur:elliptic} and biorthogonal rational functions with
self-similar spectral properties \cite{spi-zhe:spectral,spi-zhe:classical,
spi-zhe:gevp,spi:special}.
\end{remark}

\begin{definition}
The series $\sum_{n\in\mathbb{N}}c_n$ and $\sum_{n\in\mathbb{Z}}c_n$
are called {\em totally elliptic} hypergeometric series if
$h(n)=c_{n+1}/c_n$ is an elliptic function of {\em all free
parameters} entering it (except of the parameter $z$ by which one can
always multiply $h(n)$) with equal periods of double periodicity.
\end{definition}

\begin{theorem}
The most general (in the sense of a maximal number of independent
free parameters) totally elliptic theta hypergeometric series coincide with
the well-poised balanced theta hypergeometric series $_{r}E_{r-1}$
(in the unilateral case) and $_rG_r$ (in the bilateral case) for $r>2$.
Totally elliptic series are automatically modular invariant.
\end{theorem}
\begin{proof}
It is sufficient to prove this theorem for the bilateral series since
the unilateral series can be obtained afterwards by a simple reduction.
Ellipticity in $n$ leads to $h(n)$ of the form
$$
h(n)=\frac{[n+u_1,\ldots, n+u_r]}{[n+v_1,\ldots, n+v_r]}\, z,
$$
with the free parameters $u_1,\ldots,u_r$ and
$v_1,\ldots,v_{r}$ satisfying the
balancing condition $u_1+\ldots+u_r=v_1+\ldots+v_{r}$. From such a
representation it is evident that there is a freedom in the shift of
parameters by an arbitrary constant: $u_m\to u_m+u_0$, $v_m\to v_m+u_0$,
$m=1,\ldots, r$, which does not spoil the balancing condition.

Let us determine now the maximal possible number of independent
variables in the totally elliptic hypergeometric series. In general one
can denote as $a_l,\, l=1,\ldots, L,$ a set of free parameters of the
elliptic hypergeometric series in which the series is doubly periodic
with some periods.
Then $u_m=\sum_{k=1}^L \alpha_{mk}a_k+\beta_{m}$ and
$v_m=\sum_{k=1}^L \gamma_{mk}a_k+\delta_{m}$ are some linear combinations
of $a_l$ with integer
coefficients $\alpha_{mk}, \gamma_{mk}$. However, because of the possibilities to
change variables we can take a number of $u_m$ starting, say, from $u_1$
and $v_m$ as $a_1,\ldots,a_L$ and demand the double periodicity in these
parameters themselves. Since the minimal order of elliptic function is equal
to 2, the function $h(n)$ should have at least two zeros and two poles
(or one double zero or pole) in $u_1$. Double zeros and poles ask for
additional constraints, i.e. to a reduction of the number of free
parameters, and we discard such a possibility.

Let us assume that $u_r$ depends linearly on $u_1$ and suppose that
$u_1,\ldots,u_{r-1}$ are independent variables. Then it is evident that
all denominator parameters $v_m,\, m=1,\ldots,r,$ cannot contain
additional independent variables. Indeed, if it would be so, then,
inevitably, this parameter should show up at least in one $\theta$-function
in the numerator, which cannot happen by the assumption.

So, $L=r-1$ is the maximal possible number of independent variables in
the totally elliptic hypergeometric series and
$u_r$ together with $v_m,\, m=1,\dots,r,$ depend linearly on
$u_k,\, k=1,\ldots,r-1.$ Because of the permutational invariance in
the latter variables one must have $u_r=\alpha\sum_{k=1}^{r-1}u_k+\beta$,
where $\alpha, \beta$ are some numerical coefficients to be determined
(evidently, $\alpha$ must be an integer).
As to the choice of $v_m$ the unique option guaranteeing
permutational invariance of the product $\prod_{m=1}^r[x+v_m]$ in
$u_1,\ldots,u_{r-1}$ is the following one
$$
v_m=\gamma\sum_{k=1}^{r-1}u_k+\delta u_m+\rho,
\quad m=1,\ldots,r-1,
$$
and $v_r=\mu\sum_{k=1}^{r-1}u_k+\nu$,
where $\gamma, \delta, \rho, \mu, \nu$ are some numerical parameters
($\gamma, \delta, \mu$ must be integers).
This is the most general choice of $v_m$ since all other permutationally
invariant combinations of $u_m$ require products of more than $r$
theta functions. Substitution of the taken ansatz into the balancing
condition yields $1+\alpha=(r-1)\gamma+\delta+\mu$ and
$\beta=(r-1)\rho+\nu$, which guarantees invariance of $h(x)$ under the
shift $u_k\to u_k+\sigma^{-1}$ and cancels the sign factor emerging
from the shift $u_k\to u_k+\tau\sigma^{-1}$. A bit cumbersome but
technically straightforward analysis of the condition of cancellation
of the factors of the form $e^{-2\pi i\sigma u}$ yields the equations
$$
\delta^2=1,\qquad \alpha^2=(r-1)\gamma^2+2\gamma\delta+\mu^2
$$
and
$$
\alpha\beta=(r-1)\gamma\rho+\rho\delta+\mu\nu.
$$
The constraint generated by the cancellation of the factors of the form
$e^{-\pi i\tau}$ appears to be irrelevant and it will not be indicated.

Let $\delta=1$. Then two equations upon the coefficients $\alpha,\gamma,
\mu$ yield that either $\gamma=0$ or $\gamma(r-1)(r-2)/2+\mu(r-1)=1.$ Since
$\gamma$ and $\mu$ are integers, the second case cannot be valid
(the integers on the left hand side are proportional to $r-1$ whereas on
the right hand such proportionality does not take place for $r>2$).
The choice $\gamma=0$ leads to $\alpha=\mu$ and, from other two
equations, one gets $\rho=0$ and $\beta=\nu$. As a result, $h(n)=1$,
i.e. we get a trivial solution which is discarded.

Let $\delta=-1$. Solution of the taken equations gives uniquely
$$
\alpha=\frac{\gamma r}{2}, \quad
\mu=1-\frac{\gamma (r-2)}{2}, \quad
\beta=\frac{\rho r}{2},\quad \nu=-\frac{\rho (r-2)}{2},
$$
where $\gamma$ is an integer (for odd $r$ it must be an even
number) and $\rho$ is an arbitrary parameter. Arbitrariness of $\rho$
seems to contradict the statement that there are no new independent
parameters in the variable $u_r$. This paradox is resolved in the following
way. Let us first make the shift $\rho\to \rho -\gamma\sum_{k=1}^{r-1}u_k$.
It is easy to see that this leads to the removal of $\gamma$
from $h(n)$, i.e. it is a fake parameter. Denote now $\rho\equiv 2u_0$
and make the shifts $u_m\to u_m+u_0$, $m=1,\ldots,r-1$. As a result, one gets
\begin{equation}\label{vwp-theorem}
h(n)=\prod_{m=1}^{r-1}\frac{[n+u_0+u_m]}{[n+u_0-u_m]}\,
\frac{[n+u_0-\sum_{k=1}^{r-1}u_k]}{[n+u_0+\sum_{k=1}^{r-1}u_k]}\, z,
\end{equation}
i.e. the parameter $u_0$ plays the same role as $n$ and ellipticity
in it is evident. It is not difficult to recognize in (\ref{vwp-theorem})
the most general expression for $c_{n+1}/c_n$ of well-poised and balanced
theta hypergeometric series. Thus we have proved that ellipticity in
all free parameters in $h(n)$ leads uniquely to the balancing and
well-poisedness conditions.

Let us prove now that the totally elliptic hypergeometric series are
automatically modular invariant. For this is it sufficient to check
that the sum of squares of the parameters $u$ entering elliptic numbers
$[n+u]$ in the numerator of $h(n)$ (\ref{vwp-theorem}) and denominator
coincide. The numerator parameters generate the sum
$$
\sum_{k=1}^{r-1}u_k^2 +\Bigl(-\sum_{k=1}^{r-1} u_k\Bigr)^2
$$
which is trivially equal to the sum appearing from the denominator:
$$
\sum_{k=1}^{r-1}(-u_k)^2+\Bigl(\sum_{k=1}^{r-1} u_k\Bigr)^2,
$$
i.e. modular invariance is automatic. (Note that well-poised theta
hypergeometric series without balancing condition are not modular
invariant).

All the considerations given above were designed for the bilateral
$_{r}G_r$ series case, but a passage to the unilateral well-poised
and balanced $_rE_{r-1}$-series is done by a simple specification
of one of the parameters. One has to set $u_{r-1}=u_0-1$ and then
shift $u_0\to u_0+1/2,\, u_m\to u_m-1/2,\, m=1,\ldots, r-2$. This brings
$h(n)$ to the form $h(n)=z\prod_{m=0}^{r-1}[n+u_0+u_m]/[n+1+u_0-u_m]$,
where we introduced anew the $u_{r-1}$ parameter through the
relation $\sum_{m=0}^{r-1}u_m=r/2$.
\end{proof}

\begin{remark}
We have replaced $t_0$ by $t_0^2$ in the definition of very-well-poised
elliptic hypergeometric series $_{r+1}E_r$ in order to have $p$-shift
invariance in the variable $t_0$ (or ellipticity in $u_0$). Otherwise
there would be $p$-shift invariance in $t_0^{1/2}$.
\end{remark}

Thus we have found interesting origins of the balancing and well-poisedness
conditions for the plain and basic hypergeometric series calling for a
revision of these notions. However, origins of the very-well-poisedness
condition remain unknown. Probably elliptic functions $h(n)$ obeying
such a constraint have some particular arithmetic properties. An indication
on this is given by the transformation and summation formulas for some
of the very-well-poised balanced theta hypergeometric series of the
unit argument.

The theta hypergeometric series $_{r}E_s$ and $_rG_s$
are defined as formal infinite series. However, because of the
quasiperiodicity of the theta functions it is not a simple task
to determine their convergence and this problem will not be considered
here. It can be shown that for some choice of parameters the radius of
convergence $R$ of the balanced $_{r+1}E_r$-series is equal to 1.
If $R$ is the radius of convergence of the very-well-poised
$_{r+2}E_{r+1}$-series without balancing condition,
i.e. if these infinite series are well defined for $|z|<R$, then from
the representation (\ref{G/E}) it follows that  the $_rG_r$-series
converge for $\left|{q^{r-8}}/{R\prod_{m=1}^{r-4}t_m^2}\right|<|z|<R.$
A rigorous meaning to the $_rE_s$-series can be given by imposing some
truncation conditions.
The theta hypergeometric series truncates if for some $m$,
\begin{equation}\label{e-trunc1}
u_m=-N-K\sigma^{-1}-M\tau\sigma^{-1}, \qquad N\in\mathbb{N},\quad
K, M\in\mathbb{Z},
\end{equation}
or in the multiplicative form
\begin{equation}\label{e-trunc2}
t_m=q^{-N}p^{-M}, \qquad N\in\mathbb{N},\quad M\in\mathbb{Z}.
\end{equation}
The well-poised elliptic hypergeometric series are double periodic in
their parameters with the periods $\sigma^{-1}$ and $\tau\sigma^{-1}$.
Therefore, these truncated series do not depend on the integers $K, M$.

The top-level identity in the theory of basic hypergeometric series is
the four term Bailey identity for non-terminating $_{10}\Phi_9$
very-well-poised balanced series of the unit argument \cite{gas-rah:basic}.
In the terminating case there remains only two terms.
In \cite{fre-tur:elliptic} Frenkel and Turaev have proved an
elliptic generalization of the Bailey identity in the terminating case.
In our notations it looks as follows
\begin{eqnarray}\nonumber
{_{12}E_{11}}(t_0;t_1,\dots,t_7;q,p;1) &=&
\frac{\theta(qt_0^2,qs_0/s_4,qs_0/s_5,q/t_4t_5;p;q)_N}
{\theta(qs_0^2,qt_0/t_4,qt_0/t_5,q/s_4s_5;p;q)_N}   \\
&& \times {_{12}E_{11}}(s_0;s_1,\dots,s_7;q,p;1),
\label{ft-bailey}\end{eqnarray}
where it is assumed that $\prod_{m=0}^7t_m=q^2$, $t_0t_6=q^{-N}$,
$N\in\mathbb{N}$, and
$$
s_0^2=\frac{qt_0}{t_1t_2t_3},\quad s_1=\frac{s_0t_1}{t_0},\quad
s_2=\frac{s_0t_2}{t_0},\quad s_3=\frac{s_0t_3}{t_0},
$$
$$
s_4=\frac{t_0t_4}{s_0},\quad s_5=\frac{t_0t_5}{s_0},\quad
s_6=\frac{t_0t_6}{s_0},\quad s_7=\frac{t_0t_7}{s_0}.
$$

If one sets $t_2t_3=q$, then the left-hand side
of (\ref{ft-bailey}) becomes a terminating $_{10}E_9$-series, and in the
series on the right-hand side one gets $s_1=1$, i.e. only its first
term is different from zero. This gives the Frenkel-Turaev sum---an
elliptic generalization of the Jackson's sum for terminating
very-well-poised balanced $_8\Phi_7$-series \cite{gas-rah:basic}.
After diminishing the indices of $t_{4,5,6,7}$ by two
it takes the following form:
\begin{eqnarray}\nonumber
\lefteqn{ {_{10}E_9}(t_0;t_1,\dots, t_5;q,p;1) } && \\ && \makebox[4em]{}
= \frac{\theta (qt_0^2;p;q)_N\prod_{1\leq r<s\leq 3}
      \theta (q/t_rt_s;p;q)_N}{\theta (q/t_0t_1t_2t_3;p;q)_N
\prod_{r=1}^3\theta (qt_0/t_r;p;q)_N},
\label{ft-sum}\end{eqnarray}
where the parameters $t_r$ are assumed to satisfy the balancing
condition $\prod_{r=0}^5 t_r =q$ and the truncation
condition $t_0t_4=q^{-N}$, $N\in \mathbb{N}$.

\begin{remark}
Due to the clarification of the relation of very-well-poisedness
condition with the general structure of theta hypergeometric series,
starting from this paper we change notations for the elliptic
hypergeometric series in the generalizations of Bailey and
Jackson identities. The symbols $_{10}E_9$ and $_8E_7$ in the papers
\cite{spi-zhe:spectral,spi-zhe:classical,spi-zhe:gevp,spi:solitons,
spi:special,die-spi:elliptic,die-spi:selberg,die-spi:modular,
die-spi:elliptic2} read in the current notations as $_{12}E_{11}$
and $_{10}E_9$ respectively.
\end{remark}

\begin{remark}
Despite of the double periodicity, infinite totally elliptic hypergeometric
functions are not elliptic functions of $u_r$ since they have infinitely
many poles in the parallelogram of periods. Indeed, some
of the poles in $t_s,\; s=1,\ldots,r-4$, are located at $t_s=t_0q^{n+1}p^m,$
where $n\in\mathbb{N}$ and $m\in\mathbb{Z}$.
If $q^k\neq p^l$ for any $k, l\in\mathbb{N},$ then, evidently, there are
infinitely many integers $n$ and $m$ such that $t_s$ stays, say, within
the bounds $|p|< |t_s|< 1.$ This means that there are infinitely many
poles in the parameters $t_s$ in this annulus.
\end{remark}

\begin{remark}
The symbols $_rE_s$ and $_rG_s$ were chosen for denoting the one-sided
and bilateral theta hypergeometric series in order to make them as
close as possible to the standard notations $_{r}F_s$ and $_rH_s$ used for
one-sided and bilateral plain hypergeometric series respectively.
The letter $``E"$ refers also to the word ``elliptic". To the author's taste
greek symbols $_{r}\Phi_s$ and $_r\Psi_s$ used for denoting basic
hypergeometric series fit well enough into the aesthetics created
by the sequence of letters $E, \Phi, F, G, \Psi, H$.
\end{remark}

\section{Multiple elliptic hypergeometric series}

Following the definition of theta hypergeometric series
of a single variable one can consider formal multiple sums of
quasiperiodic combinations of Jacobi $\theta_1$-functions depending
of more than one summation index. However, we shall limit ourselves
only to the multiple elliptic hypergeometric series.

\begin{definition}
The formal series
$$
\sum_{\lambda_1,\ldots,\lambda_n=0}^\infty c(\lambda_1,\ldots,\lambda_n)
\quad \text{or}
\sum_{\lambda_1,\ldots,\lambda_n=-\infty}^\infty
c(\lambda_1,\ldots,\lambda_n),
$$
and
$$
\sum_{\lambda_1,\ldots,\lambda_n=0
\atop \lambda_1\leq \ldots\leq \lambda_n}^\infty
c(\lambda_1,\ldots,\lambda_n)
\quad \text{or}
\sum_{\lambda_1,\ldots,\lambda_n=-\infty
\atop \lambda_1\leq \ldots\leq \lambda_n}^\infty
c(\lambda_1,\ldots,\lambda_n)
$$
are called {\em multiple elliptic hypergeometric series} if:
a) the coefficients $c(\mathbf{\lambda})$ are symmetric with respect
to an action of permutation group $\mathcal{S}_n$ upon the
summation variables $\lambda_1,\ldots,\lambda_n$ and the free
parameters entering $c(\mathbf{\lambda})$;
b) for all $k=1,\ldots,n$ the functions
$$
h_k(\mathbf{\lambda})=\frac{
c(\lambda_1,\ldots,\lambda_k+1,\ldots,\lambda_n)
}{c(\lambda_1,\ldots,\lambda_k,\ldots,\lambda_n)},
$$
are elliptic in $\lambda_k,\, k=1,\ldots,n,$
considered as complex variables. These series are called totally elliptic
if, in addition, the functions $h_k(\mathbf{\lambda})$ are
elliptic in all free parameters except of the free multiplication factors.
\end{definition}

Suppose that $h_k(\mathbf{\lambda})$ is symmetric in $\lambda_1,\ldots,
\lambda_{k-1},\lambda_{k+1},\ldots,\lambda_n$. Then, using the results
of the one variable analysis, it is not difficult to see that the most
general expression for the coefficients $c(\mathbf{\lambda})$ is:
\begin{equation}\label{c-general}
c(\mathbf{\lambda})=
\prod_{k=1}^n\Biggl(\prod_{1\leq i_1<\ldots<i_k\leq n} \prod_{m=1}^{r_k}
\frac{[u_{km}]_{\lambda_{i_1}+\ldots+\lambda_{i_k} } }
{[v_{km}]_{\lambda_{i_1}+\ldots+\lambda_{i_k} } }\Biggr)
z_1^{\lambda_1}\cdots z_n^{\lambda_n},
\end{equation}
where
$$
\sum_{k=1}^n C_{n-1}^{k-1} \sum_{m=1}^{r_k} (u_{km}-v_{km})=0.
$$
However, if the action of $\mathcal{S}_n$ permutes $\lambda_1,\ldots,
\lambda_n$ and simultaneously free parameters entering
$c(\mathbf{\lambda})$ other than $z_1,\ldots,z_n$, then the situation is
richer, e.g. more general combinations of $\lambda_k$ are allowed than it
is indicated in (\ref{c-general}). We shall not go into further analysis
of general situation but pass to some particular examples.

Currently there are two known examples of multiple elliptic hypergeometric
series leading to some constructive identities (multivariable analogues of
the Frenkel-Turaev $_{10}E_9$-summation formula). The first one corresponds
to an elliptic extension of the Aomoto-Ito-Macdonald type of series
\cite{aom:elliptic,ito:theta,mac:constant}.
Its structure is read off from the following multivariable generalization
of the Frenkel-Turaev summation formula considered in \cite{war:summation,
die-spi:elliptic,ros:proof}.
Let $N\in\mathbb{N}$ and the parameters $t, t_r\in\mathbb{C},\,
r=0,\ldots,5,$ are constrained by the balancing condition
$t^{2n-2}\prod_{r=0}^5t_r=q$ and the truncation condition
$t^{n-1}t_0t_4=q^{-N}.$
Then one has the following theta-functions identity
\begin{eqnarray}\nonumber
&&
\sum_{0\leq \lambda_1\leq \lambda_2\leq \cdots \leq \lambda_n\leq N}
q^{\sum_{j=1}^n\lambda_j} t^{2 \sum_{j=1}^n (n-j)\lambda_j}
\\  \nonumber && \makebox[8em]{}\times
\prod_{1\leq j<k\leq n} \Biggl(
\frac{\theta(\tau_k\tau_jq^{\lambda_k+\lambda_j},
\tau_k\tau_j^{-1}q^{\lambda_k-\lambda_j};p)}
{\theta(\tau_k\tau_j,\tau_k\tau_j^{-1};p)} \\  \nonumber
&& \makebox[8em]{}\times
\frac{\theta(t\tau_k\tau_j;p;q)_{\lambda_k+\lambda_j}}
     {\theta(qt^{-1}\tau_k\tau_j;p;q)_{\lambda_k+\lambda_j}}
\frac{\theta(t\tau_k\tau_j^{-1};p;q)_{\lambda_k-\lambda_j}}
     {\theta(qt^{-1}\tau_k\tau_j^{-1};p;q)_{\lambda_k-\lambda_j}}
\Biggr) \\  \nonumber
&&\makebox[3em]{} \times \prod_{j=1}^n
\Biggl( \frac{\theta(\tau_j^2q^{2\lambda_j};p)}{\theta(\tau_j^2;p)}
\prod_{r=0}^5 \frac{\theta(t_r\tau_j;p;q)_{\lambda_j}}
     {\theta(qt_r^{-1}\tau_j;p;q)_{\lambda_j}}\Biggr)  \\
&&\makebox[3em]{}
=\frac{\theta(qt^{n+j-2}t_0^2;p;q)_N
      \prod_{1\leq r <s \leq 3} \theta(qt^{1-j} t_r^{-1}t_s^{-1};p;q)_N}
      {\theta(qt^{2-n-j}\prod_{r=0}^3t_r^{-1};p;q)_N
      \prod_{r=1}^3 \theta(qt^{j-1}t_0t_r^{-1};p;q)_N}.
\label{multi-1}\end{eqnarray}
Here the parameters $\tau_j$ are related to $t_0$ and $t$ as follows:
$\tau_j=t_0t^{j-1}$, $j=1,\ldots ,n$.
Note that the series coefficients $c(\mathbf{\lambda})$ are symmetric
with respect to the simultaneous permutation of the variables
$\lambda_j$ and $\lambda_k$ and the parameters $\tau_j$ and $\tau_k$
for arbitrary $j\neq k$ (for this one has to assume that $\tau_j$ are
independent variables).

\begin{theorem}
The series standing on the left hand side of (\ref{multi-1}) is a totally
elliptic multiple hypergeometric series.
\end{theorem}
\begin{proof}
Ratios of the series coefficients yield
\begin{eqnarray}\nonumber
h_l(\mathbf{\lambda})=\prod_{j=1}^{l-1}
\frac{\theta(\tau_j\tau_lq^{\lambda_j+\lambda_l+1},\tau_j^{-1}\tau_l
q^{\lambda_l+1-\lambda_j}, t\tau_j\tau_lq^{\lambda_j+\lambda_l},
t\tau_j^{-1}\tau_lq^{\lambda_l-\lambda_j};p)}
{\theta(\tau_j\tau_lq^{\lambda_j+\lambda_l},\tau_j^{-1}\tau_l
q^{\lambda_l-\lambda_j}, t^{-1}\tau_j\tau_lq^{\lambda_j+\lambda_l+1},
t^{-1}\tau_j^{-1}\tau_lq^{\lambda_l+1-\lambda_j};p)}
\\  \nonumber
\times \prod_{k=l+1}^n\frac{\theta(\tau_k\tau_lq^{\lambda_k+\lambda_l+1},
\tau_k\tau_l^{-1}q^{\lambda_k-\lambda_l-1},t\tau_k\tau_lq^{\lambda_k+
\lambda_l},t^{-1}\tau_k\tau_l^{-1}q^{\lambda_k-\lambda_l};p)}
{\theta(\tau_k\tau_lq^{\lambda_k+\lambda_l},\tau_k\tau_l^{-1}
q^{\lambda_k-\lambda_l}, t^{-1}\tau_k\tau_lq^{\lambda_k+\lambda_l+1},
t\tau_k\tau_l^{-1} q^{\lambda_k-\lambda_l-1};p)}
\\
\times qt^{2(n-l)}
\frac{\theta(\tau_l^2q^{2\lambda_l+2};p)}{\theta(\tau_l^2;p)}
\prod_{m=0}^5\frac{\theta(t_m\tau_lq^{\lambda_l};p)}
{\theta(t_m^{-1}\tau_lq^{\lambda_l+1};p)}.
\label{h-AIM}\end{eqnarray}

Using the equalities (\ref{fun-rel}) one can easily check the ellipticity
of this $h_l(\mathbf{\lambda})$ in $\lambda_i$ for $i< l$ and $i>l$
(for this it is simply necessary to see that $h_l$ does not change after
the replacement of $q^{\lambda_i}$ by $pq^{\lambda_i}$). For the
ellipticity in $\lambda_l$ itself one has an essentially longer computation.
The replacement of $q^{\lambda_l}$ by $pq^{\lambda_l}$ in the
product $\prod_{j=1}^{l-1}$ yields a multiplier $t^{-4(l-1)}$.
The product $\prod_{k=l+1}^n$ yields the multiplier $t^{-4(n-l)}$. The
remaining part of $h_l$ generates the factor $q^{-4}\prod_{m=0}^5
qt_m^{-2}$. The product of all these three factors takes the form
$q^2t^{-4(n-1)}\prod_{m=0}^5 t_m^{-2}$ and it is equal to 1 due to the
balancing condition.

So, we have found that the taken series is indeed a multiple elliptic
hypergeometric series. Let us prove now its total ellipticity or
$p$-shift invariance in the parameters $t_m,\, m=0,\ldots,4$ and $t$.
The $p$-shift invariance in the parameters $t_1,\ldots,t_4$ follows
from the balancing condition in the same way as in the single variable
series case. Consider the $t_0\to pt_0$ shift.
The product $\prod_{j=1}^{l-1}$ yields a multiplier $t^{-4(l-1)}$
and the product $\prod_{k=l+1}^n$ yields the factor $t^{-4(n-l)}$. The
remaining part of $h_l$ generates the factor $q^2\prod_{m=0}^5
t_m^{-2}$. The product of all these three multipliers is equal to 1.

Finally, the shift $t\to pt$ calls for the most complicated computation.
The product $\prod_{j=1}^{l-1}$ yields a complicated multiplier
$(q^{2\lambda_l+1}t^{2(l-1)}\tau_l^2p^{2l-3})^{2(1-l)}$.
The product $\prod_{k=l+1}^n$ generates no less complicated
expression $(q^{2\lambda_l+1}t^{2(l-1)}\tau_l^2p^{2(l-1)})^{2(l-n)}$.
The remaining part of $h_l$ leads to the following factor (after the use
of the balancing condition): $(q^{2\lambda_l+1}t^{2(l-1)}\tau_l^2)^{2(n-1)}
p^{(4n-6)(l-1)}$. The product of all these three multipliers yields 1.
Thus we have proved the total ellipticity of the taken type of series.
\end{proof}

The second example of multiple series corresponds to an elliptic
generalization of the Milne-Gustafson type multiple basic hypergeometric
series \cite{mil:multidimensional,gus:macdonald,den-gus:q-beta},
which are, in turn, $q$-analogues of the Hollman, Biedenharn, and Louck
plain multiple hypergeometric series \cite{hbl:hypergeometric}.
Its structure is read off from the following summation formula suggested
in \cite{die-spi:modular}.
Let $q^n\neq p^m$ for $n,m\in\mathbb{N}$. Then
for the parameters $t_0,\ldots ,t_{2n+3}$ subject to the balancing condition
$q^{-1}\prod_{r=0}^{2n+3}t_r=1$ and the truncation conditions
$q^{N_j}t_jt_{n+j}=1,\, j=1,\ldots ,n,$ where $N_j\in\mathbb{N},$
one has the identity
\begin{eqnarray}\nonumber
&& \sum_{\stackrel{0\leq \lambda_j\leq N_j}{j=1,\ldots ,n}}
q^{\sum_{j=1}^n j\lambda_j}
\prod_{1\leq j<k\leq n}
\frac{\theta (t_jt_kq^{\lambda_j+\lambda_k},
              t_jt_k^{-1}q^{\lambda_j-\lambda_k} ;p)}
     {\theta (t_jt_k,t_jt_k^{-1} ;p)}  \\
\label{multi-2}
&& \quad\qquad\qquad\qquad\times \prod_{1\leq j\leq n} \Biggl(
\frac{\theta (t_j^2q^{2\lambda_j};p)}{\theta (t_j^2;p)}
\prod_{0\leq r\leq 2n+3} \frac{\theta (t_jt_r;p;q)_{\lambda_j}}
     {\theta (qt_jt_r^{-1};p;q)_{\lambda_j}}\Biggr)  \\
&&= \theta (qa^{-1}b^{-1},qa^{-1}c^{-1},qb^{-1}c^{-1};p;q)_{N_1+\cdots +N_n}
\nonumber\\
&&\qquad\times\prod_{1\leq j<k\leq n}
\frac{\theta (qt_jt_k;p;q)_{N_j}\theta (qt_jt_k ;p;q)_{N_k}}
     {\theta (qt_jt_k;p;q)_{N_j+N_k}} \nonumber\\
&&\qquad \times \prod_{1\leq j\leq n}
\frac{\theta (qt_j^2;p;q)_{N_j}}
     {\theta (qt_ja^{-1},qt_jb^{-1},
 qt_jc^{-1},q^{1+N_1+\cdots+N_n-N_j}t_j^{-1}a^{-1}b^{-1}c^{-1};p;q)_{N_j}},
\nonumber\end{eqnarray}
where $a\equiv t_{2n+1}$, $b\equiv t_{2n+2}$, $c\equiv t_{2n+3}$.
Note that this series coefficients $c(\mathbf{\lambda})$ are symmetric
with respect to simultaneous permutation of the variables
$\lambda_j$ and $\lambda_k$ together with the parameters $t_j$ and $t_k$
for arbitrary $j,k=1,\ldots,n,\, j\neq k$.

\begin{theorem}
The series standing on the left hand side of (\ref{multi-2}) is a
totally elliptic hypergeometric series.
\end{theorem}
\begin{proof}
Ratios of the successive series coefficients yield
\begin{eqnarray}\nonumber
h_l(\mathbf{\lambda})= \prod_{j=1}^{l-1}
\frac{\theta(t_jt_lq^{\lambda_j+\lambda_l+1},t_jt_l^{-1}
q^{\lambda_j-\lambda_l-1};p)}{\theta(t_jt_lq^{\lambda_j+\lambda_l},
t_jt_l^{-1}q^{\lambda_j-\lambda_l};p)}  \\  \nonumber
\times
\prod_{k=l+1}^{n}
\frac{\theta(t_lt_kq^{\lambda_l+\lambda_k+1},t_lt_k^{-1}
q^{\lambda_l+1-\lambda_k};p)}{\theta(t_lt_kq^{\lambda_l+\lambda_k},
t_lt_k^{-1}q^{\lambda_l-\lambda_k};p)}  \\  \times
q^l\frac{\theta(t_l^2q^{2\lambda_l+2};p)}
{\theta(t_l^2q^{2\lambda_l};p)} \prod_{m=0}^{2n+3}
\frac{\theta(t_lt_mq^{\lambda_l};p)}{\theta(t_lt_m^{-1}q^{\lambda_l+1};p)}.
\label{h-MG}\end{eqnarray}
It is easy to check the ellipticity of this expression in $\lambda_j$
for $j<l$ and $j>l$. Ellipticity in $\lambda_l$ itself follows from
a more complicated computation. Namely, the change of $q^{\lambda_j}$
to $pq^{\lambda_j}$ leads to additional multipliers in
$h_l(\mathbf{\lambda})$: $q^{2(1-l)}$---from the product
$\prod_{j=1}^{l-1}$, $q^{2(l-n)}$---from the product $\prod_{k=l+1}^n$, and
$q^{-4}\prod_{m=0}^{2n+3}qt_m^{-2}$---from the rest of
$h_l(\mathbf{\lambda})$. Multiplication of these three
expressions gives 1 due to the balancing condition.

Ellipticity in the parameters $t_m,\, m=0,\ldots,2n+3,$ is checked
separately for $m<l, m>l$ and $m=l$. The first two cases are easy
enough and do not worth of special consideration. The replacement
$t_l\to pt_l$ leads to the following multipliers:
$q^{2(1-l)}$---from the product $\prod_{j=1}^{l-1}$,
$q^{2(l-n)}$---from the product $\prod_{k=l+1}^n$, and
$q^{-2n}\prod_{m=0}^{2n+3}t_m^{-2}$---from the rest of
$h_l(\mathbf{\lambda})$.
The balancing condition guarantees again that the total multiplier is
equal to 1. Thus we have proved total ellipticity of this series as well.
\end{proof}

\begin{remark}
Modular invariance of the series (\ref{multi-1}) and (\ref{multi-2})
has been established in \cite{die-spi:elliptic} and \cite{die-spi:modular}
respectively. We conjecture that all totally elliptic multiple
hypergeometric series are automatically modular invariant similar to
the one-variable series situation. One can introduce general notions of
well-poised and very-well-poised multiple theta hypergeometric series,
with the given above examples being counted as very-well-poised series,
but we shall not discuss this topic in the present paper.
\end{remark}

 The author is indebted to J.F. van Diejen and A.S. Zhedanov for
 a collaboration in the work on elliptic hypergeometric series
 and for useful discussions of this paper. Valuable comments and
 encouragement from G.E. Andrews, R. Askey, and A. Berkovich
 are highly appreciated as well.

\bibliographystyle{amsplain}

\end{document}